\documentclass[12pt,reqno, a4paper]{amsart}

\usepackage{amsmath,amsthm,amsfonts,amssymb,amsxtra,appendix,bookmark,dsfont,bm,color}
\usepackage[shortlabels]{enumitem}
\usepackage{graphicx}
\usepackage[utf8]{inputenc}

\usepackage{tcolorbox}
\theoremstyle{plain}
\newtheorem{theorem}{Theorem}

\theoremstyle{remark}
\newtheorem{remark}{Remark}
\theoremstyle{definition}

\newcounter{step} 

\usepackage{etoolbox}
\AtBeginEnvironment{proof}{\setcounter{step}{0}}

\def\bq{\begin{eqnarray}}
\def\eq{\end{eqnarray}}
\def\bqq{\begin{align*}}
\def\eqq{\end{align*}}

\def\eps{\varepsilon}

\def\R {\mathbb{R}}

\def\N {\mathcal{N}}

\def\R {\mathbb{R}}

\def\N {\mathbb{N}}

\def\d{{\, \rm d}}

\def\Cr{C_{{\rm rad}, 2^*}}
\def\Ci{C_{{\rm rad},\infty}}
\usepackage{pdfpages}
\usepackage[normalem]{ulem}

\title[Hardy--Sobolev interpolation inequalities]{Hardy--Sobolev interpolation inequalities}

\author[C. Dietze]{Charlotte Dietze}
\address[C. Dietze]{Department of Mathematics, LMU Munich, Theresienstrasse 39, 80333 Munich, Germany} 
\email{dietze@math.lmu.de}

\author[P.T. Nam]{Phan Th\`anh Nam}
\address[P.T. Nam]{Department of Mathematics, LMU Munich, Theresienstrasse 39, 80333 Munich, Germany} 
\email{nam@math.lmu.de}

\begin{document}

\begin{abstract} We derive a family of interpolation estimates which improve Hardy's inequality and cover the Sobolev critical exponent. We also determine all optimizers among radial functions in the endpoint case and discuss open questions on nonrestricted optimizers.  
\end{abstract}

\maketitle

\section{Introduction}

The classical Hardy inequality states that for every dimension $d\ge 3$, 
\bq \label{eq:Hardy}
h[u]:=\int_{\R^d} |\nabla u(x)|^2 \d x - \frac{(d-2)^2}{4} \int_{\R^d}\frac{|u(x)|^2}{|x|^2} \d x \ge 0, \quad \forall u\in \dot{H}^1(\R^d). 
\eq
In this short article we are interested in interpolation inequalities involving the quadratic form $h[u]$ and $L^p$-norms of $u$. A classical result in this direction is the Gagliardo--Nirenberg type inequality
\bq \label{eq:BV1997}
h[u]^{\theta} \|u\|_{L^2}^{2(1-\theta)}  \ge C \|u\|_{L^q}^2, \quad \forall u\in H^1(\R^d)
\eq
for a constant $C>0$ independent of $u$, which holds for every  
$$d\ge 3, \quad 2<q<2^*=\frac{2d}{d-2}, \quad \theta=d\left(  \frac 1 2  - \frac 1 q\right).$$
The inequality \eqref{eq:BV1997} can be deduced from the results of Brezis and V\'azquez \cite[Theorem 4.1 and Extension 4.3]{BV1997}; see \cite{VZ2000} for related results. The bound \eqref{eq:BV1997} can be also derived from Sobolev's embedding theorem and the kinetic estimate 
\bq \label{eq:F2009}
h[u]^{\theta} \|u\|_{2}^{2(1-\theta)}  \ge C \|(-\Delta)^{s/2}u\|_{2}^2, \quad \forall u\in H^1(\R^d),
\eq
 for $s\in(0,1)$ and $\theta=\theta(s)$, which was proved by Frank \cite[Theorem 1.2]{F2009}. Both of \eqref{eq:BV1997} and \eqref{eq:F2009} have been extended to the fractional Laplacian in \cite{F2009}, motivated by applications in the asymptotic behavior of large Coulomb systems \cite{SSS2010} and the stability of  relativistic matter  \cite{FLS2008,F2009}. 

Note that the restriction $q<2^*$ in \eqref{eq:BV1997} is necessary, namely the quadratic form $h[u]$ is really weaker than $\|\nabla u\|_{L^2}^2$. Here we are interested in a replacement of \eqref{eq:BV1997} which covers the critical power $q=2^*$, with the expense that the $L^2$-norm is replaced by the energy associated with the inverse square potential. We have

\begin{theorem}[Hardy-Sobolev interpolation inequality]\label{thm:1} If $d=3$ and $\theta=1/3$, then the inequality
\bq \label{eq:Hardy-improved-1}
\left(  \int_{\R^d} |\nabla u|^2 - \frac{(d-2)^2}{4} \sup_{y\in \R^d} \int_{\R^d}\frac{|u(x)|^2}{|x-y|^2} \d x  \right)^{\theta}  \left( \sup_{y\in \R^d} \int_{\R^d}\frac{|u(x)|^2}{|x-y|^2} \d x  \right)^{1-\theta}  \ge C \| u\|_{L^{2^*}}^2
\eq
holds with a constant  $C=C(d,\theta)>0$  independent of $u\in \dot H^1(\R^d)$. Moreover, \eqref{eq:Hardy-improved-1} does not hold if $d\ge 4$ or if 
$\theta \ne 1/3$. 
\end{theorem}

\begin{remark} The bound \eqref{eq:Hardy-improved-1} is invariant under translations and dilations. Note that for the first term on the left-hand side, Hardy's inequality \eqref{eq:Hardy} is equivalent to  
$$
\int_{\R^d} |\nabla u|^2 - \frac{(d-2)^2}{4}\sup_{y\in \R^d} \int_{\R^d}\frac{|u(x)|^2}{|x-y|^2} \d x \ge 0, \quad \forall u\in \dot{H}^1(\R^d). 
$$
For the second term, it is important to include $\sup_{y\in \R^d}$ since otherwise this term can be made arbitrarily small by translation $u\mapsto u(\cdot - z)$ with $|z|\to \infty$. 
\end{remark} 

\begin{remark} For all $d\ge 3$ and $1-2/d \le \theta \le 1$ we have
\bq \label{eq:Hardy-improved-1'}
\left(  \int_{\R^d} |\nabla u|^2 \right)^{\theta}  \left( \sup_{y\in \R^d} \int_{\R^d}\frac{|u(x)|^2}{|x-y|^2} \d x  \right)^{1-\theta}  \ge C \| u\|_{L^{2^*}}^2,\quad \forall u\in H^1(\R^d).
\eq
This is a consequence of the improved Sobolev inequality involving Morrey norms
\bq \label{eq:Hardy-improved-1''}
\left(  \int_{\R^d} |\nabla u|^2 \right)^{\theta}  \left( \sup_{R>0,x\in \R^d} R^{-2} \int_{B(x,R)} |u|^2 \right)^{1-\theta}  \ge C \| u\|_{L^{2^*}}^2,\quad \forall u\in \dot H^1(\R^d),
\eq
which was proved by Palatucci--Pisante \cite[Theorem 1]{PP2014}, using subtle weighted $L^p$-estimates for Riesz potentials in \cite{SW1992} and Calder\'on-Zygmund type techniques in the spirit of the Fefferman--Phong argument \cite{F1983}. The bound \eqref{eq:Hardy-improved-1''} is helpful to obtain the compactness of minimizing sequences of the critical Sobolev inequality; see \cite[Theorem 3]{PP2014} for details. In contrast, our inequality \eqref{eq:Hardy-improved-1} is stronger than \eqref{eq:Hardy-improved-1'} and it only holds for the special case $d=1/\theta=3$. 

\end{remark} 

In the next result, we extend  \eqref{eq:Hardy-improved-1} by replacing the gradient term $\| \nabla u\|_{L^2}$ by $\|\nabla u\|_{L^p}$, as well as replacing  the $L^{2^*}$-norm by the $L^{p^*,r}$-Lorentz norm. Recall that (see \cite[Definition 1.4.6 and Proposition 1.4.9]{G-book})
$$
\| u \|_{L^{p,r}}=\| u \|_{{p,r}}= \begin{cases}
\left( p\int_0^\infty s^{r-1} |\{|u|>s\}|^{r/p} \d s \right)^{1/r}, \quad 0<r<\infty, \\
\sup_{s>0} s |\{ |u|>s\}|^{1/p} ,\quad r=\infty.
\end{cases}
$$

\begin{theorem} [Hardy-Sobolev inequalities with Lorentz norms]  \label{thm:2} 
Let $d\ge 2$, $p\in [2,d)$, $p^*= pd/(d-p)$, $r\in [p, \infty]$ and
\begin{align} \label{eq:range-theta}
\theta\in\left[\frac{p}{\min(r,p^*)},\frac{1}{p}-\frac{1}{r}\right].
\end{align} 
Then 
\begin{align} \label{eq:general-p-r}
\left( \int_{\R^d} |\nabla u|^p - \left( \frac{d-p}{p}\right)^p \sup_{y\in \R^d} \int_{\R^d} \frac{|u(x)|^p}{|x-y|^p} \d x \right)^{\theta}\left( \sup_{y\in \R^d} \int_{\R^d} \frac{|u(x)|^p}{|x-y|^p} \d x \right)^{1-\theta}  \ge C \|u\|_{L^{p^*,r}}^p
\end{align}
 with a constant $C=C(d,p,r,\theta)>0$ independent of $u\in \dot{W}^{1,p}(\R^d)$. The bound \eqref{eq:general-p-r} does not hold if  $\theta<p/\min(r,p^*)$ (with arbitrary $p\ge 2$), or if $\theta>1/p-1/r$ and $p=2$. 
In particular, when $p=2$, the range of $\theta$ in \eqref{eq:range-theta} is optimal. 
\end{theorem}

\medskip

Theorems \ref{thm:1} and \ref{thm:2} naturally lead  to the question of determining optimizers of the relevant inequalities. We expect that in the non-endpoint cases 
$$\frac{p}{\min(r,p^*)}< \theta < \frac{1}{p}-\frac{1}{r},$$
the existence of optimizers of \eqref{eq:general-p-r} follows from the standard concentration compactness method. Below we focus on the endpoint cases. While the existence of optimizers in this case is open in general, we are able to give a partial answer under the restriction to radial functions. We will limit ourselves to the choice $p=2$ and $r\in \{2^*,\infty\}$, for which the right-hand side of \eqref{eq:general-p-r} becomes either the usual $L^{2^*}$-norm or the $L^{2^*,\infty}$-weak norm.  The relevant functional space is 
\begin{align*}
\dot{H}^1_{\rm rad}(\R^d)&=\{u \in \dot H^1(\R^d): u \text { is radially symmetric}\}.
\end{align*}
In the radial case, we can work directly with the quadratic form $h[u]$ in \eqref{eq:Hardy}. We have

\begin{theorem}[Radial optimizers]\label{thm:3} Let $d\ge 3$ and $p=2$.  

\begin{enumerate}[(i)]
    \item Let $r=2^*=2d/(d-2)$ and $\theta=1/p-1/r=1/d$. Then all optimizers of the inequality 
 \bq \label{eq:Hardy-radial}
    h[u] ^\theta  \left(  \int_{\R^d}\frac{|u(x)|^2}{|x|^2} \d x  \right)^{1-\theta} \ge C_{\rm rad,2^*} \| u\|_{L^{2^*}}^{2}, \quad \forall u\in \dot{H}^1_{\rm rad}(\R^d)
\eq
are given by the family  
\begin{align} \label{eq:Terracini}
u_\eta(x) = \frac{1} { (|x|^{1-\eta} (1+|x|^{2\eta}))^{(d-2)/2}}, \quad \eta \in (0,\infty),
\end{align} 
up to dilation $u_\eta(x) \mapsto a u_\eta (bx)$ with $a\in \mathbb{C}, b>0.$ Furthermore, \eqref{eq:Hardy-radial}  does not hold if $\theta \ne 1/d$. 

    \item Let $r=\infty$ and $\theta=1/p-1/r=1/2$. Then all optimizers of the inequality
        \bq \label{eq:Hardy-radial-inf}
        h[u]^\theta  \left(  \int_{\R^d}\frac{|u(x)|^2}{|x|^2} \d x  \right)^{1-\theta} \ge C_{\rm rad,\infty} \| u\|_{L^{2^*,\infty}}^{2},  \quad \forall u\in \dot{H}^1_{\rm rad}(\R^d)
\eq
are given by the family 
\begin{align} \label{eq:u-f-putback-intro}
u_c(x)=\begin{cases}
|x|^{c-d/2+1}\,\, &,\, |x|\le 1,\\
|x|^{-c-d/2+1}\,\,&, \, |x| > 1, 
\end{cases} \quad c\in (0, d/2-1],
\end{align}
up to dilation.  
\end{enumerate}
\end{theorem}

\begin{remark}
The classification of all optimizers in Theorem \ref{thm:3} (i) is consistent with Terracini's study in \cite{Terracini1996} where all radial positive solution of the Euler-Lagrange equation 
\begin{align} \label{eq:HS-equation}
-\Delta u(x) - \frac{(d-2)^2}{4} (1-\eta^2) \frac{u(x)}{|x|^2}  = u^{2^*-1}(x), \quad x\in \R^d,
\end{align} 
with a given constant $\eta>0$, were derived. In particular, according to \cite[Eq. (4.6)]{Terracini1996}, the only regular solutions (i.e., belonging to $L^{2^*}$), up to rescaling, are of the form $(d(d-2) \eta^2)^{(d-2)/4} u_\eta$ with $u_\eta$ given in \eqref{eq:Terracini}. 
\end{remark}

\begin{remark} We leave the following {\bf open questions}: Do optimizers of  \eqref{eq:Hardy-improved-1} exist? And if exist, are they radial? The same questions for the simpler inequality \eqref{eq:Hardy-improved-1'} in the endpoint case $\theta=1-2/d$ remain unsolved. Note that both in \eqref{eq:Hardy-improved-1} and  \eqref{eq:Hardy-improved-1'}, all quantities scale in the same way. Moreover, by taking several bubbles travelling far from each other, one can construct optimizing sequences that do not converge weakly to a nonzero limit after any choice of dilations and translations. It seems that a novel concentration-compactness argument will be needed to resolve the existence problem of optimizers for these inequalities.

\end{remark}

We prove Theorem \ref{thm:3} in Section \ref{sec:radial}, and then prove Theorems \ref{thm:1}  and \ref{thm:2} in Section \ref{sec:general}.

\subsection*{Acknowledgements.} Part of the work has been done when CD was a visiting reseacher at the Institut des Hautes Études and she would like to thank Laure Saint-Raymond for her support and hospitality. We are grateful to Rupert L. Frank for helpful suggestions, leading to our consideration of Lorentz norms in Theorem \ref{thm:2}.  We also thank Hajer Bahouri, Morris Brooks, Albert Cohen, Patrick Gérard, Kihyun Kim, Mathieu Lewin,  Hoai-Minh Nguyen, Julien Sabin, Nikita Simonov, Thomas S\o rensen, Jakob Stern, Hanne van den Bosch and Jean van Schaftingen for interesting and inspiring discussions. 
We acknowledge the support from the Deutsche Forschungsgemeinschaft through the DFG project Nr. 426365943. CD also acknowledges the support from the Jean-Paul Gimon Fund and from the Erasmus+ programme.

\section{Radial case} \label{sec:radial}

In this section we prove Theorem \ref{thm:3}. Let $u\in H^1_{\rm rad}(\R^d)$ with  $d\ge 3$. 

\medskip
\noindent
{\bf Proof of (i):} Using the ground state representation for Hardy's inequality (see e.g. \cite[Eq. (2.14)]{FS08}), we denote
\begin{align} \label{eq:gs-u-f}
u(x) = \frac{f(|x|)}{|x|^{(d-2)/2}}
\end{align}
and rewrite 
\begin{align*}
\int_{\R^d} |u|^{2^*} & =  |\mathbb{S}^{d-1} | \int_0^\infty \frac{ |f(r)|^{2^*}}{r} \d r,\quad \int_{\R^d} \frac{|u(x)|^2}{|x|^2} \d x  =|\mathbb{S}^{d-1} | \int_0^\infty \frac{ |f(r)|^2}{r} \d r,\\
\int_{\R^d} |\nabla u|^{2}  &- \frac{(d-2)^2}{4}\int_{\R^d} \frac{|u(x)|^2}{|x|^2} \d x= \int_{\R^d} \frac{|\nabla f(x)|^2 }{|x|^{d-2}} \d x =  |\mathbb{S}^{d-1} |\int_{0}^\infty |f'(r)|^2 r \d r.
\end{align*}
Here $|\mathbb{S}^{d-1} |$ is the surface area of the unit sphere in $\R^d$. Thus \eqref{eq:Hardy-radial} is  equivalent to 
\begin{align}\label{eq:1D-form-proof-0}
\left( \int_0^\infty r |f'(r)|^2 {\rm d} r \right)^{\theta} \left( \int_0^\infty \frac{|f(r)|^2}{r} {\rm d} r \right)^{1-\theta} \ge \Cr |\mathbb{S}^{d-1} |^{2/2^*-1} \left( \int_0^\infty \frac{|f(r)|^{2^*}}{r} \d r \right)^{2/2^*}.   
\end{align}

The bound \eqref{eq:1D-form-proof-0} can be interpreted as a Caffarelli--Kohn--Nirenberg type inequality \cite{CKN1984}, namely
\begin{align}\label{eq:1D-form-proof}
\||x|^{1/2} \nabla f \|_{L^2(\R_+)}^{\theta} \| |x|^{-1/2} f\|_{L^2(\R_+)}^{1-\theta} \ge \sqrt{\Cr} |\mathbb{S}^{d-1} |^{1/2^*-1/2} \| |x|^{\gamma} f \|_{L^{r}(\R_+)}  
\end{align}
with $r=2^*=2d/(d-2)$, $\gamma=-1/r$. Actually it is a limiting case as $1/r+\gamma/n=0$  in dimension $n=1$, which does not seem available from the literature (see \cite[Theorem 3.1]{NS2018}, \cite[Theorem 1.2]{NS2019} and \cite[Theorem 2.2]{MN2021} for recent results in the limiting case). 

Inspired by \cite{Terracini1996}, we make the following changes of variables
$$
f(r) =\psi(\log r), \quad s=\log r \in \R, \quad \d s=\frac{\d r}{r},
$$
which give 
\begin{align*}
\int_0^\infty \frac{|f(r)|^2}{r} \d r &= \int_{\R} |\psi(s)|^2 \d s , \quad
\int_0^\infty \frac{|f(r)|^{2^*}}{r} \d r = \int_{\R} |\psi(s)|^{2^*} \d s ,\\
\int_0^\infty |f'(r)|^2 r \d r & = \int_{\R} |\psi'(s)|^2 \d s. 
\end{align*}
Therefore, \eqref{eq:1D-form-proof-0} is equivalent to the Gagliardo--Nirenberg interpolation inequality
\begin{align}\label{eq:LT1D-gen}
\left( \int_{\R} |\psi'(s)|^2 \d s \right)^{\theta} \left( \int_{\R} |\psi(s)|^2 \d s \right)^{1-\theta} \ge \Cr |\mathbb{S}^{d-1} |^{2/2^*-1} \left( \int_{\R} |\psi(s)|^{2^*} \d s  \right)^{2/2^*}.   
\end{align}
The optimal constant of the one-dimensional inequality \eqref{eq:LT1D-gen} was already obtained by Nagy in 1941 \cite{Nagy1941}, with $2^*=2d/(d-2)$ replaced by a general positive power. The existence and uniqueness of optmizers of the analogue of \eqref{eq:LT1D-gen} in higher dimensions are also well-known; we refer to the classical works of Weinstein \cite{Weinstein1983} and Kwong \cite{Kwong} for instance. The uniqueness of optimizers of \eqref{eq:LT1D-gen} can be translated straightforwardly to the classification of optmizers of  \eqref{eq:Hardy-radial}  as stated in Theorem \ref{thm:3} (i); we refer to \cite[Eq. (4.6)]{Terracini1996} for a similar analysis. 

\begin{remark}
In the special case $d=3$ (which is relevant to Theorem \ref{thm:1}), the interpolation inequality \eqref{eq:LT1D-gen} with $2^*=6$ goes back to the (1D, one-body) Lieb--Thirring inequality \cite{LT-76} as well as Keller's lower bound on the lowest eigenvalue of the Schr\"odinger operator  $-\d^2/\d x^2 + V(x)$ on $L^2(\R)$ \cite{K-61}; see also \cite[Section 2]{BL-04} for a simple derivation of the optimal constant in this special case. 
\end{remark}

\medskip
\noindent
{\bf Unique choice of $\theta$ for \eqref{eq:Hardy-radial}:} Consider \eqref{eq:1D-form-proof-0} with the trial function
\begin{align} \label{eq:example-f-reps}
f(r)=\begin{cases}
r^{\eps},\quad r\in (0,1],\\
r^{-\eps},\quad r\in [1,\infty),
\end{cases}
\end{align}
where $\eps>0$ is a parameter. Then we have
\begin{align}
\int_0^\infty r |f'(r)|^2 \d r = \eps,\quad \int_0^\infty \frac{|f(r)|^2}{r} \d r =\frac{1}{\eps},\quad \int_0^\infty \frac{|f(r)|^{2^*}}{r} \d r =\frac{2}{2^* \eps}.
\end{align}
Therefore, \eqref{eq:1D-form-proof-0} requires to have 
\begin{align} \label{eq:eps-eps-eps}
\eps^{\theta} \eps^{-(1-\theta)} \ge C  \eps^{- 2/2^*}
\end{align}
for all $\eps>0$. By letting $\eps\to0$ and $\eps\to\infty$, we find that 
$$
\theta= \frac{1-2/2^*}{2}= \frac{1}{d}. 
$$
Thus, \eqref{eq:Hardy-radial} holds only if $\theta=1/d$.

\medskip
\noindent
{\bf Proof of (ii):} Let us consider \eqref{eq:Hardy-radial-inf}. Using again the ground state representation \eqref{eq:gs-u-f},  we have
\begin{align} \label{eq:u-2*-inf}
\|u\|_{L^{2^*,\infty}} &= \sup_{t>0} t |\{x: |u(x)| >t \}|^{1/2^*} = \sup_{t>0} t |\{x: |f(x)| >t |x|^{d/2^*} \}|^{1/2^*} \nonumber \\
&\le \sup_{t>0} t |\{x: \|f\|_{L^\infty} >t |x|^{d/2^*} \}|^{1/2^*} = |B(0,1)|^{1/2^*} \|f\|_{L^\infty}.
\end{align}
Here $|B(0,1)|$ is the volume of the unit ball in $\R^d$. Therefore, \eqref{eq:Hardy-radial-inf} holds if we can show that 
\begin{equation} \label{eq:eq:Hardy-radial-inf-v2}
    \left( \int_0^\infty r |f'(r)|^2 \d r \right)  \left( \int_0^\infty \frac{|f(r)|^2}{r} \d r \right)  \ge  \Ci^2|\mathbb{S}^{d-1} |^{-2}|B(0,1)|^{4/2^*}\|f\|_{L^\infty}^{4}.
\end{equation}
From \eqref{eq:gs-u-f} and $u\in \dot H^1(\R^d)$, we deduce that the function $f:(0,\infty)\to [0,\infty)$ is continuous and satisfies 
$$\lim_{r\to 0}f(r)=0,\quad \lim_{r\to \infty}f(r)=0.$$
Therefore, up to dilation, we may without loss of generality assume that 
$$f(1)=\|f\|_{L^\infty}.$$
By the Cauchy-Schwarz inequality and the fundamental theorem of calculus, we have
\begin{align}\begin{split}\label{eq:cs-inf-f}
   &\quad \left( \int_0^\infty r |f'(r)|^2 \d r \right)  \left( \int_0^\infty \frac{|f(r)|^2}{r} \d r \right)  \ge \left( \int_0^\infty |f'(r) f(r)|   \d r \right)^2\\
   &
   = \left(\frac{1}{2} \left(\int_0^1 \left\vert\frac{\partial}{\partial r} (|f(r)|^2)\right\vert   \d r + \int_1^\infty \left\vert\frac{\partial}{\partial r} (|f(r)|^2)\right\vert   \d r  \right)\right)^2
   \ge |f(1)|^4 = \|f\|_{L^\infty}^{4}.  
\end{split}
\end{align}
Thus \eqref{eq:eq:Hardy-radial-inf-v2} holds, and consequently \eqref{eq:Hardy-radial-inf} holds, with 
\begin{align} \label{eq:opt-Ci}
\Ci^2|\mathbb{S}^{d-1} |^{-2}|B(0,1)|^{4/2^*}=1. 
\end{align}

To have the equality in \eqref{eq:Hardy-radial-inf}, we need to ensure all equalities in \eqref{eq:u-2*-inf} and \eqref{eq:cs-inf-f}. The bound \eqref{eq:cs-inf-f} contains two inequalities where the first equality occurs if there exists a constant $c> 0$ such that 
$$r |f'(r)|^2=c\frac{|f(r)|^2}{r},\quad \text{\rm a.e. } r\in (0,\infty),$$
while the second equality occurs if $|f|^2$ is monotone increasing on $(0,1)$ and monotone decreasing on $(1,\infty)$. Thus, we have all equalities in \eqref{eq:cs-inf-f} if and only if
\begin{align} \label{eq:f-rc}
f(r)=\begin{cases}
r^{c}\quad &,\, r\in (0,1],\\
r^{-c}\quad &,\, r\in [1,\infty). 
\end{cases}
\end{align}

It remains to determine the range of $c$ in \eqref{eq:f-rc} to get the equality in \eqref{eq:u-2*-inf}. If $0<c \le (d-2)/2$, then the function 
\begin{align} \label{eq:u-f-putback}
u(x)= f(|x|) |x|^{-\frac{d-2}{2}}= \begin{cases}
|x|^{c-\frac{d-2}{2}}\quad &,\, |x|\le 1,\\
|x|^{-c-\frac{d-2}{2}}\quad &, \, |x|\ge 1 
\end{cases}
\end{align}
is radially symmetric decreasing, and the equality in \eqref{eq:u-2*-inf} occurs since 
\begin{align} \label{eq:u-2*-inf-lower}
\|u\|_{L^{2^*,\infty}} &= \sup_{t>0} t |\{x: |u(x)| >t \}|^{1/2^*} \ge |\{x: |u(x)| > 1 \}|^{1/2^*} = |B(0,1)|^{1/2^*}.
\end{align}
On the other hand, if $c>(d-2)/2$, then the inequality in \eqref{eq:u-2*-inf} is strict: since $u$ defined in \eqref{eq:u-f-putback} is bounded by $1$, we have 
\begin{align*}
\|u\|_{L^{2^*,\infty}}^{2^*} &= \sup_{0<t<1} t^{2^*} |\{x: |u(x)| >t \}| \nonumber\\
&=  \sup_{0<t<1} t^{2^*} \left( |\{|x|\le 1: |u(x)| >t \}| + |\{|x|>1: |u(x)| >t \}|  \right) \nonumber\\
&= \sup_{0<t<1} t^{2^*} \left( |\{ x: 1 \ge |x|> t^{\frac 2 {2c-(d-2)}} \}| + |\{x: t^{ -\frac 2 {2c+d-2} } > |x| >1\}|  \right) \nonumber\\
&= |B(0,1)|  \sup_{0<t<1} t^{2^*} \left( t^{ -\frac {2d} {2c+d-2} }  - t^{\frac {2d} {2c-(d-2)}}  \right) < |B(0,1)|,
\end{align*}
where the latter estimate can be easily seen using the fact that
$$
2^* = \frac{2d}{d-2} >  \frac {2d} {2c+d-2}.
$$

Thus, in summary, \eqref{eq:Hardy-radial-inf} holds with the optimal constant $\Ci$ given in \eqref{eq:opt-Ci}, and all optimizers are uniquely characterized up to dilation by \eqref{eq:u-f-putback} with $c\in (0, (d-2)/2]$. 

The proof of Theorem \ref{thm:3} is complete. 
\qed

\bigskip

\section{General case} \label{sec:general}

In this section we prove Theorem \ref{thm:1} and Theorem \ref{thm:2}. 

\subsection{Proof of Theorem \ref{thm:1}}
We divide the proof into two parts. First, we prove \eqref{eq:Hardy-improved-1} for $d=3=1/\theta$. Then we show that the condition  $d=3=1/\theta$ is necessary. 

\begin{proof}[Proof of \eqref{eq:Hardy-improved-1} for $d=3=1/\theta$] Let $u\in \dot{H}^1(\R^d)$ and denote
$$
A = \int_{\R^d} |\nabla u|^2 , \quad B= \frac{(d-2)^2}{4} \sup_{y\in \R^d} \int_{\R^d}\frac{|u(x)|^2}{|x-y|^2} \d x. 
$$
We consider two cases.

\medskip
\noindent
{\bf Case 1:} $(1-\theta)A \ge B$. Then using 
$$
B \ge \frac{(d-2)^2}{4} \sup_{r>0,y\in \R^d} \int_{B(y,r)}\frac{|u(x)|^2}{|x-y|^2} \d x \gtrsim \sup_{r>0, y \in \R^d} \frac{1}{r^2} \int_{B(y,r)} |u|^2,
$$
and $A-B\ge \theta A$,  we conclude from \eqref{eq:Hardy-improved-1''} (see \cite[Theorem 1]{PP2014}) that
$$
(A-B)^{\theta} B^{1-\theta} \gtrsim \|\nabla u\|_{L^2}^\theta \left( \sup_{r>0, y \in \R^d} \frac{1}{r^2} \int_{B(y,r)} |u|^2 \right)^{1-\theta}  \gtrsim  \|u\|_{L^{2^*}}^2.
$$

\medskip
\noindent
{\bf Case 2:} $(1-\theta)A <B$.  Then $B\mapsto (A-B)^\theta B^{1-\theta}$ is monotone decreasing since
$$
\frac{\d }{\d B} ( (A-B)^\theta B^{1-\theta} )= ((1-\theta) A-B) (A-B)^{\theta-1} B^{-\theta}<0. 
$$
Moreover, by  the Hardy-Littlewood  and P\'olya-Szeg\"o   rearrangement inequalities (see e.g. \cite[Lemma 1.6 and Theorem 4.7]{Burchard2009}) we have
$$
B^* = \frac{(d-2)^2}{4}\int_{\R^d} \frac{|u^*(x)|^2}{|x|^2} \d x \ge B, \quad A^*= \|\nabla u^*\|_{L^2}^2 \le A,
$$
where $u^*$ denotes the radially symmetric decreasing rearrangement of $u$. Therefore,   
$$
(A-B)^\theta B^{\theta} \ge (A-B^*)^\theta (B^*)^{1-\theta} \ge (A^*-B^*)^\theta (B^*)^{1-\theta}. 
$$
Thus it remains to consider \eqref{eq:Hardy-improved-1} in the case when $u$ is radially symmetric decreasing. In this case, 
\begin{equation} \label{eq:sup-rev}
    \sup_{y\in \R^d} \int_{\R^d}\frac{|u(x)|^2}{|x-y|^2} \d x=\int_{\R^d}\frac{|u(x)|^2}{|x|^2}. 
\end{equation}
by the Hardy-Littlewood rearrangement inequality. Now the desired bound has been already proved in Theorem \ref{thm:3} (i).

Thus in all cases, \eqref{eq:Hardy-improved-1} holds for $d=3=1/\theta$. 
\end{proof}

\begin{proof}[Proof of the necessity of $d=3=1/\theta$] Let us show that \eqref{eq:Hardy-improved-1} fails if $d\ge 4$ or if $\theta\ne 1/3$. First, the necessity of $\theta\le 1/d$ can be seen from the radial case as explained in Theorem \ref{thm:3} (i). To be precise, we consider the example in \eqref{eq:example-f-reps} with $\eps>0$ small. In this case, $u$ is radially symmetric decreasing, and hence \eqref{eq:sup-rev} holds.
Therefore, \eqref{eq:Hardy-improved-1} requires \eqref{eq:eps-eps-eps} for $\eps>0$ small, which implies that $\theta\le 1/d$.

\medskip

In order to complete the proof, we consider another (non-radial) example. Fix $\varphi\in C_c^\infty\backslash\{0\}$, $z\in \mathbb{R}^d\backslash\{0\}$ and choose 
$$u_N(x)= \sum_{n=1}^N \varphi (x+ n N z)$$
with $N\to \infty$. Then by replacing $u$ by $u_N$, we find that 
$$
A \sim N, \quad B \sim 1, \quad \|u\|_{L^{2^*}}^{2^*} \sim N. 
$$
for large $N$. Therefore, \eqref{eq:Hardy-improved-1} requires 
$$
N^{\theta} \gtrsim N^{2/2^*}, 
$$
which implies that $
\theta   \ge 2/2^* = 1-2/d.
$
Combining with the upper bound $\theta\le 1/d$ and the constraint $d\ge 3$, we find that the only possibility is $d=3$ and $\theta=1/3$. \end{proof}

The proof of Theorem \ref{thm:1} is complete.

\subsection{Proof of Theorem \ref{thm:2}} \label{sec:thm3} We first prove \eqref{eq:general-p-r} and then explain the necessity of the constraint of $\theta$. 

\begin{proof}[Proof of \eqref{eq:general-p-r}]
Let $d\ge 2$ and $p\in [2,d)$. Assume that $r\in [p, \infty]$ and 
\begin{align} \label{eq:theta-int}
\frac{p}{\min(r,p^*)} \le \theta \le \frac{1}{p}-\frac{1}{r}.
\end{align}
Fix a small constant $\eps=\eps(d,p,r,\theta)\in(0,1)$. We consider two cases.

\medskip
\noindent
{\bf Case 1:} Assume
\begin{equation}\label{eq:case1-eps}
(1-\eps)\int_{\R^d} |\nabla u|^p \ge  \left( \frac{d-p}{p}\right)^p \sup_{y\in \R^d} \int_{\R^d} \frac{|u(x)|^p}{|x-y|^p} \d x.
\end{equation}
From \cite[Theorem 1]{PP2014}, see also \cite[Eq. (1.4)]{Sch}, we have
\begin{align} \label{eq:Sch}
\left( \int_{\R^d} |\nabla u|^p \right)^{p/p^*} \left( \sup_{y\in \R^d, R>0} \frac{1}{R^p} \int_{B(y,R)} |u(x)|^p \d x \right)^{1-p/p^*}  \gtrsim \|u\|_{L^{p^*}}^p.
\end{align}
A simplified proof of \eqref{eq:Sch} based on sharp maximal functions can be obtained by following the analysis in \cite{TN}.  We can extend this bound to the Lorentz norm on $L^{p^*,r}$ with $r\in [p,p^*]$, namely 
\begin{align} \label{lem:C}
\left( \int_{\R^d} |\nabla u|^p \right)^{p/r} \left( \sup_{y\in \R^d, R>0} \frac{1}{R^p} \int_{B(y,R)} |u(x)|^p \d x \right)^{1-p/r}  \gtrsim \|u\|_{L^{p^*,r}} ^p. 
\end{align}
To prove \eqref{lem:C}, let us use the standard dyadic decomposition: recalling that $\varphi:\R\to \R$ is a smooth function supported on the annulus $1/2\le |t| \le 2$ such that 
$$
\sum_{j\in \mathbb{Z}} \varphi(2^{-j}t)=1,\quad \forall t\in \R\backslash\{0\},
$$
we write 
$$
u =\sum_{j\in \mathbb{Z}} u_j, \quad u_j = u \varphi(2^{-j} |u|). 
$$
Then 
\begin{align*}
\|u\|_{L^{p^*,r}}^{p^*} & \sim  \left(  \sum_{j\in \mathbb{Z}} \|u_j\|_{L^{p*}}^r \right)^{p^*/r}  \nonumber\\
& \lesssim \left(  \sum_{j\in \mathbb{Z}} \left( \int_{\R^d} |\nabla u_j|^p \right)^{r/p^*} \left( \sup_{y\in \R^d, R>0} \frac{1}{R^p} \int_{B(y,R)} |u_j(x)|^p \d x \right)^{r/d} 
 \right)^{p^*/r} ,
\end{align*}
where we used \eqref{eq:Sch} for $u_j$ and the fact that $(1-p/p^*)r/p=r/d$ (as $1/p-1/d=1/p^*$). On the other hand, using 
$$
\sup_{y\in \R^d, R>0} \frac{1}{R^p} \int_{B(y,R)} |u_j(x)|^p \d x \lesssim \min\left\{ \int_{\R^d} |\nabla u_j|^p,  \sup_{y\in \R^d, R>0} \frac{1}{R^p} \int_{B(y,R)} |u(x)|^p \d x \right\}
$$
and splitting the power 
$$r/d= r(1/p-1/p^*) = (1- r/p^*) + (r/p-1),$$
we find that
\begin{align*}
&\left( \sup_{y\in \R^d, R>0} \frac{1}{R^p} \int_{B(y,R)} |u_j(x)|^p \d x \right)^{r/d} \nonumber\\
&\lesssim \left( \int_{\R^d} |\nabla u_j|^p \right)^{1-r/p^*} \left( \sup_{y\in \R^d, R>0} \frac{1}{R^p} \int_{B(y,R)} |u(x)|^p \d x \right)^{r/p-1}.   
\end{align*}
Here we used the constraint $r\in [p,p^*]$ to ensure that both $(1- r/p^*)$ and $(r/p-1)$ are nonnegative. Thus, we obtain 
\begin{align*}
\|u\|_{L^{p^*,r}}^{p^*}  &\lesssim \left(  \sum_{j\in \mathbb{Z}} \left( \int_{\R^d} |\nabla u_j|^p \right) \left( \sup_{y\in \R^d, R>0} \frac{1}{R^p} \int_{B(y,R)} |u(x)|^p \d x \right)^{r/p-1} 
 \right)^{p^*/r} \\
 &\lesssim  \left(  \left( \int_{\R^d} |\nabla u|^p \right) \left( \sup_{y\in \R^d, R>0} \frac{1}{R^p} \int_{B(y,R)} |u(x)|^p \d x \right)^{r/p-1} 
 \right)^{p^*/r} ,
\end{align*}
which is equivalent to \eqref{lem:C}. 

From \eqref{eq:Sch}, \eqref{lem:C} and the obvious bound $\|u\|_{L^{p^*}}^p \gtrsim \|u\|_{p^*, r}^p$ for $r\ge p^*$, we get
\begin{equation}
    \left( \int_{\R^d} |\nabla u|^p  \right)^{\frac{p}{\min(r,p^*)}} \left( \sup_{y\in \R^d} \int_{\R^d} \frac{|u(x)|^p}{|x-y|^p} \d x \right)^{1-\frac{p}{\min(r,p^*)}} \gtrsim \|u\|_{p^*, r}^p.
\end{equation}
By Hardy's inequality and the condition $\theta\ge p/\min(r,p^*)$, we also get
\begin{equation}\label{eq:theta-Morrey}
    \left( \int_{\R^d} |\nabla u|^p  \right)^{\theta} \left( \sup_{y\in \R^d} \int_{\R^d} \frac{|u(x)|^p}{|x-y|^p} \d x \right)^{1-\theta} \gtrsim \|u\|_{p^*, r}^p.
\end{equation}
Combining \eqref{eq:case1-eps} and \eqref{eq:theta-Morrey}, we obtain
\begin{align}
\begin{split}
&\left( \int_{\R^d} |\nabla u|^p - \left( \frac{d-p}{p}\right)^p \sup_{y\in \R^d} \int_{\R^d} \frac{|u(x)|^p}{|x-y|^p} \d x \right)^{\theta} \left( \sup_{y\in \R^d} \int_{\R^d} \frac{|u(x)|^p}{|x-y|^p} \d x \right)^{1-\theta}  \\
&\ge \left(\eps \int_{\R^d} |\nabla u|^p  \right)^{\theta} \left( \sup_{y\in \R^d} \int_{\R^d} \frac{|u(x)|^p}{|x-y|^p} \d x \right)^{1-\theta}   \gtrsim  \|u\|_{p^*, r}^p. 
\end{split}
\end{align}

\medskip
\noindent
{\bf Case 2:} Assume 
$$
(1-\eps)\int_{\R^d} |\nabla u|^p \le  \left( \frac{d-p}{p}\right)^p \sup_{y\in \R^d} \int_{\R^d} \frac{|u(x)|^p}{|x-y|^p} \d x.
$$
Then by a monotonicity argument as in the proof of Theorem \ref{thm:1}, we can reduce to the case when $u$ is radially symmetric decreasing. 
 We have for radially symmetric decreasing functions $u$, 
\begin{equation}\label{eq:lor-rad-dec-f}
\|u\|_{q,s}^s\sim \int_{\R^d} \frac{|u(x)|^s}{|x|^\alpha} \d x
\end{equation}
if $s/q=1-\alpha/d$ (see e.g. \cite[Lemma 4.3]{FS08} for the case $q=p^*$, $s=p$). In particular, for $q=p^*=dp/(d-p)$ and $s=\alpha=d$, we obtain 
\begin{equation}\label{eq:lor-rad-dec-f-formula}
\|u\|_{p^*,s}^s \sim  \int_{\R^d} \frac{|f(x)|^s}{|x|^d} \d x \sim \int_0^\infty \frac{|f(r)|^s}{r} \d r
\end{equation}
with $u(x)=|x|^{1-d/p} f(x)$. Moreover,  
$$
\int_{\R^d} \frac{|u(x)|^p}{|x|^p} \d x = \int_{\R^d}  \frac{|f(x)|^p}{|x|^d} \d x \sim \int_0^\infty \frac{|f(r)|^p}{r} \d r, 
$$
and by the ground state representation \cite[Eq. (2.14)]{FS08} (here we use that $p\ge2$)
$$
\int_{\R^d} |\nabla u|^p - \left( \frac{d-p}{p}\right)^p \int_{\R^d} \frac{|u(x)|^p}{|x|^p} \d x  \gtrsim 
\int_{\R^d} \frac{|\nabla f(x)|^p}{|x|^{d-p}} \d x \sim \int_0^\infty |f'(s)|^p s^{p-1} \d s.
$$
By H\"older's inequality, we have, with $1/p+1/p'=1$, 
\begin{align}
&\left( \int_0^\infty |f'(s)|^p s^{p-1} \d s \right)^{1/p} \left(  \int_0^\infty \frac{|f(s)|^p}{s} \d s \right)^{1/p'} \nonumber  \\
&\ge \int_0^\infty |f'(s)| |f(s)|^{p-1} = p^{-1} \int_0^\infty | (f^p(s))'| \d s \gtrsim_p  \|f\|_{\infty}^p. 
\end{align}
Under the constraint \eqref{eq:theta-int}, there exists $\tilde r\in [p,r]$ such that $\theta= 1/p - 1/\tilde r$. Let $\beta\ge 0$ such that $\tilde r = p(1+\beta)$. Then
\begin{align}
&\left( \int_0^\infty |f'(s)|^p s^{p-1} \d s \right)^{\beta/p} \left(  \int_0^\infty \frac{|f(s)|^p}{s} \d s \right)^{1+\beta/p'}  \nonumber\\
&\ge p^{-1} \|f\|_{\infty}^{\beta p}  \int_0^\infty \frac{|f(s)|^p}{s} \d s  = p^{-1}\int_0^\infty \frac{|f(s)|^{\tilde r}}{s} \d s\sim \|u\|_{p^*,\tilde r}^{\tilde r}\gtrsim  \|u\|_{p^*,r}^{\tilde r},  
\end{align}
where we used \eqref{eq:lor-rad-dec-f-formula} with $s=\tilde r$ and $\tilde r\le r$. This implies \eqref{eq:general-p-r} by the choice of $\tilde r$. \end{proof}

\begin{proof}[Proof of the necessity of the range of $\theta$] Let us explain why \eqref{eq:general-p-r} fails for certain values of $\theta$ as indicated in Theorem \ref{thm:2}. First we consider  $\theta<p/\min (r,p^*)$ for general $p\ge 2$ where we split into two cases $r\ge p^*$ and $r<p^*$, and then we focus on the case $p=2$. 

\medskip

\noindent
{\bf Counterexample for the case $r\ge p^*$ and $\theta<p/\min (r,p^*)=p/p^*$. }
The idea is to consider $N\in\N$ identical bubbles that travel away from each other and to let $N\to\infty$. Let $0\ne \varphi\in C_c^\infty(B(0,1))$, $0\ne z\in \R^d$ and choose 
\begin{equation}
   u_N(x)= \sum_{n=1}^N \varphi (x+ n N z) 
\end{equation}
for every $N\in \mathbb{N}$. The translation by $Nz$ ensures that the functions $\{\varphi (\cdot + n N z)\}_{n=1}^N$ have disjoint support for $N$ large. We have 
\begin{equation}
    \sup_{y\in \R^d} \int_{\R^d} \frac{|u_N(x)|^p}{|x-y|^p} \d x =\sup_{y\in \R^d} \int_{\R^d} \frac{|\varphi(x)|^p}{|x-y|^p} \d x +o(1)\quad\text{as}\quad N\to\infty\,
\end{equation}
and
\begin{equation}
    \int_{\R^d} |\nabla u_N|^p - \left( \frac{d-p}{p}\right)^p \sup_{y\in \R^d} \int_{\R^d} \frac{|u_N(x)|^p}{|x-y|^p} \d x=N\int_{\R^d} |\nabla \varphi|^p+O(1)\quad\text{as}\quad N\to\infty\,.
\end{equation}
Moreover, it is straightforward to see that for $N$ large 
\begin{equation}
    \|u_N\|_{p^*,r}=p^* \int_0^\infty s^{r-1} |\{|u_N|>s\}|^{r/p^*} \d s = N^{r/p^*}p^* \int_0^\infty s^{r-1} |\{|\varphi|>s\}|^{r/p^*} \d s \,,
\end{equation}
if $r<\infty$, and 
\begin{equation}
    \|u_N\|_{p^*,\infty}= \sup_{s>0} s |\{ |u_N|>s\}|^{1/p^*}=N^{1/p^*}\sup_{s>0} s |\{ |\varphi|>s\}|^{1/p^*} \,.
\end{equation}
if $r=\infty$. Hence, if \eqref{eq:general-p-r} holds, then by taking $u=u_N$, we get 
\begin{equation}
    N^\theta\gtrsim N^{\frac{p}{p^*}}
\end{equation}
for $N$ large, which implies that $\theta\ge p/p^*$. 

\bigskip

\noindent
{\bf Counterexample for the case $r<p^*$ and $\theta<p/\min (r,p^*)=p/r$. }  Let $0\ne \varphi\in C_c^\infty(B(0,1))$, $0\ne z\in \R^d$ and define 
\begin{align} \label{eq:def-vN}
v_N (x) =\sum_{j=1}^N 2^j \varphi (2^{p^*j/d} (x+ jN z)) 
\end{align}
for every $N\in \mathbb{N}$. The scaling is chosen such that $\| 2^j \varphi (2^{p^*j/d} \cdot )\|_{L^{p^*}(\R^d)}=\|\varphi\|_{L^{p^*} (\R^d)}$ for all $j$, and similarly  all relevant terms in \eqref{eq:general-p-r} are invariant when changing $\varphi\mapsto 2^j \varphi (2^{p^*j/d} \cdot )$. Again, the translation by $Nz$ ensures that the supports of the functions $\{2^j \varphi (2^{p^*j/d} (\cdot + jNz) ) \}_{j=1}^N$ 
are far away from each other for $N$ large. Then 
$$
\int_{\R^d} |\nabla v_N|^p = N \int_{\R^d} |\nabla \varphi|^p, \quad \sup_{y\in \R^d} \int_{\R^d} \frac{|v_N(x)|^p}{|x-y|^p} {\rm d}x \sim \sup_{y\in \R^d} \int_{\R^d} \frac{|\varphi(x)|^p}{|x-y|^p} {\rm d}x \sim 1
$$
and
$$
\|v_N\|_{L^{p^*,r}}^{p} \sim \left( \sum_{j}  \| 2^j \varphi (2^{p^*j/d} \cdot )\|_{L^{p^*}(\R^d)}^r\right)^{p/r} \sim N^{p/r}. 
$$
Thus inserting $v_N$ in \eqref{eq:general-p-r}, we find that for $N$ large,
$$
N^{\theta} \gtrsim N^{p/r},
$$
which requires $\theta\ge p/r$. 

\bigskip

\noindent
{\bf Counterexample for the case $p=2$ and $\theta> 1/p-1/r$.} Define $u:\R^d\to\R$ by
\begin{equation}
    u(x) =  \frac{f(|x|)}{|x|^{(d-p)/p}}.
\end{equation}
where $f(r)$ is chosen as in \eqref{eq:example-f-reps} with $\eps>0$ small. Then $u$ is radially symmetric decreasing and hence \eqref{eq:sup-rev} holds. Therefore, 
\begin{align}\label{eq:upxp-comp-ex-eps}
\begin{split}
    \sup_{y\in \R^d} \int_{\R^d} \frac{|u(x)|^p}{|x-y|^p} \d x &=\int_{\R^d} \frac{|f(|x|)|^p}{|x|^d} \d x =|\mathbb{S}^{d-1} | \int_0^\infty \frac{ |f(r)|^p}{r} \d r=|\mathbb{S}^{d-1} | \frac{2}{\eps p}.
    \end{split}
\end{align}
We also have by the ground state representation for $p=2$, 
\begin{align} \label{eq:GS-repr-p=2}
     \qquad \int_{\R^d} |\nabla u|^p - \left( \frac{d-p}{p}\right)^p \sup_{y\in \R^d} \int_{\R^d} \frac{|u(x)|^p}{|x-y|^p} \d x= \int_{\R^d} \frac{|f(|x|)|^p}{|x|^{d-p}} \d x =\eps^{p-1}|\mathbb{S}^{d-1}|\frac{2}{ p} .
\end{align}
Note that the analogue of \eqref{eq:GS-repr-p=2} is more complicated for $p\ne 2$ (see \cite{FS08}), which is why our counterexample only works for $p=2$. 

By \eqref{eq:lor-rad-dec-f}, we have for $r<\infty$
\begin{align*}
    \|u\|_{L^{p^*,r}}^r&=\int_{\R^d} \frac{|f(x)|^{r}}{|x|^d} \d x = |\mathbb{S}^{d-1} | \int_0^\infty \frac{|f(s)|^{r}}{s} \d s= |\mathbb{S}^{d-1} | \frac{2}{\eps r}. 
\end{align*}
Moreover, for $r=\infty$, by \eqref{eq:u-2*-inf} and \eqref{eq:u-2*-inf-lower},  
\begin{align*}
\|u\|_{p^*,\infty} =|B(0,1)|^{1/p^*}.
\end{align*}
Therefore, if \eqref{eq:general-p-r} holds, then 
\begin{equation}
    \eps^{(p-1)\theta}\eps^{-(1-\theta)}\gtrsim \eps^{-p/r}\,,
\end{equation}
for $\eps>0$ small, which requires that 
$
    p\theta-1\ge-p/r
$ namely 
$$
    \theta\ge\frac{1}{p}-\frac{1}{r}\,.
$$

The proof of Theorem \ref{thm:2} is complete. 
\end{proof}

\medskip



\end{document}